\documentclass[smallextended,numbook,runningheads]{svjour3}     
\smartqed  
\usepackage{graphicx}
\usepackage{mathptmx}      
\usepackage{amssymb,amsmath}
\usepackage{multienum,subfloat,multirow}
\usepackage{pstricks,pst-node,pst-tree}
\psset{treemode=U,levelsep=1,treesep=0.5}

\newtheorem{theo}{Theorem}
\newtheorem{lem}[theo]{Lemma}

\newtheorem{defi}{Definition}

\newcommand{\R}{\mathbb{R}}
\newcommand{\Ih}{I^h}

\newcommand{\Ids}{1\kern-0.25em{\rm l}}

\providecommand{\E}{\operatorname{E}}
\providecommand{\N}{\mathbb{N}}
\newcommand{\bO}{\ensuremath{\mathcal{O}}}
\newcommand{\Tadd}{T^{add}}
\newcommand{\PhiH}{\Phi_H}
\newcommand{\A}{A}
\newcommand{\J}{J}
\journalname{BIT}
\begin{document}

\title{Runge-Kutta methods for third order weak
approximation of SDEs with multidimensional additive noise
}
\titlerunning{SRK methods for third order weak approximation of SDEs with additive noise}        

\author{Kristian Debrabant}


\institute{Kristian Debrabant \at
              Technische Universit\"{a}t Darmstadt, Fachbereich Mathematik, Schlo{\ss}gartenstr.7,
64289 Darmstadt, Germany \\
              \email{debrabant@mathematik.tu-darmstadt.de}\\
             \emph{Present address:}Katholieke Universiteit Leuven, Scientific Computing Research Group, Celestijnenlaan 200A, 3001 Leuven (Heverlee), Belgium
}

\date{Received: date / Accepted: date}
\maketitle

\begin{abstract}
A new class of third order Runge-Kutta methods for stochastic differential equations with additive noise is introduced. In contrast to Platen's method, which to the knowledge of the author has been up to now the only known third order Runge-Kutta scheme for weak approximation, the new class of methods affords less random variable evaluations and is also applicable to SDEs with multidimensional noise. Order conditions up to order three are calculated and coefficients of a four stage third order method are given. This method has deterministic order four and minimized error constants, and needs in addition less function evaluations than the method of Platen. Applied to some examples, the new method is compared numerically with Platen's method and some well known second order methods and yields very promising results.
\end{abstract}

\section{Introduction}
\label{aco:intro}
In many applications, e.\,g., in epidemiology and financial
mathematics, taking stochastic effects into account when modelling dynamical systems often leads to stochastic differential equations (SDEs). An important subclass of these are SDEs with additive noise in the form
\begin{equation}\label{aco:eq:SDE}
X(t)=x_0+\int_{t_0}^tg_0(s,X(s))~ds+\sum_{i=1}^mg_i(W_i(t)-W_i(t_0)).
\end{equation}
\begin{sloppypar}Here, $W(t)$ is an m-dimensional Wiener process defined on a probability space $(\Omega,\A,\mathcal{P})$,
the Borel-measurable drift $g_0:\R^d\to\R^d$ is assumed to be sufficiently differentiable and to satisfy a Lipschitz and a linear growth condition, and $g_i\in\R^d$, $i=1,\dots,m$. Then the Existence and Uniqueness Theorem \cite{karatzas91bma} applies. Examples of such systems arising in experimental psychology, turbulent diffusion, radio-astronomy and  blood clotting dynamics can be found in \cite{kloeden99nso}.
\end{sloppypar}
In recent years, the development of numerical
methods for the approximation of SDEs has become a field of
increasing interest, see e.\,g. \cite{kloeden99nso,milstein95nio} and references therein.
 Whereas strong approximation methods are
designed to obtain good pathwise solutions, see e.\,g. \cite{burrage04nmf}, weak approximation focuses on the expectation of functionals of the solution:

Let $C_P^l(\R^d, \R)$ denote the space of all $g \in
C^l(\R^d,\R)$ fulfilling a polynomial growth condition \cite{kloeden99nso}. Further, let ${\Ih} = \{t_0, t_1, \ldots, t_N\}$ with $t_0 < t_1 < \ldots < t_N =T$ be a discretization of the time interval $I=[t_0,T]$ with step sizes $h_n = t_{n+1}-t_n$ for $n=0,1, \ldots, N-1$.
\begin{defi}[weak convergence]
    A time discrete approximation
    $Y^h=(Y^h(t))_{t \in {\Ih}}$ converges weakly with order $p$ to $X$ as
    $h \rightarrow 0$ at time $t \in {\Ih}$ if for each $f \in
    C_P^{2(p+1)}(\R^d, \R)$ there exist a constant $C_f$
    and a finite $\delta_0 > 0$ such that
    \begin{equation*}
        | \E(f(Y^h(t))) - \E(f(X(t))) | \leq C_f \, h^p
    \end{equation*}
    holds for each $h \in \, ]0,\delta_0[\,$.
\end{defi}

Many approximation schemes for SDEs fall into the class of stochastic Runge-Kutta (SRK) methods. Second order SRK methods for the weak approximation of SDEs were proposed by Kloeden and Platen \cite{kloeden99nso}, Komori \cite{komori07wso}, Mackevicius and Navikas \cite{mackevicius01sow}, Tocino and Vigo-Aguiar \cite{tocino02wso}, R\"{o}{\ss}ler \cite{roessler07sor,roessler09sor}, and Debrabant and R\"{o}{\ss}ler
 \cite{debrabant08cwa,debrabant09ddi,debrabant09foe}.
 An explicit third order weak SRK method for autonomous SDEs with additive scalar noise as well as its generalization to general scalar noise have been given in Kloeden and Platen \cite{kloeden99nso}. However, the authors state there that "it remains an open and challenging task to derive simpler derivative free order 3.0 weak schemes, at least for important classes of stochastic differential equations." The present article solves this problem in the case of additive noise and overcomes also the restriction to scalar additive noise.

To do so, we consider the following class of $s$-stage SRK methods,
\begin{subequations}\label{aco:eq:SRK}
\begin{align}
Y_{n+1}&=Y_n+h_n\sum_{i=1}^s\alpha_ig_0(t_n+c_ih_n,H_i)+\sqrt{h_n}\sum_{l=1}^mg_l\J_l,\\
H_i&=Y_n+h_n\sum_{j=1}^sa_{ij}g_0(t_n+c_jh_n,H_j)+\sqrt{h_n}\sum_{l=1}^mg_l(b_{1,i}\J_l+b_{2,i}\J_{m+l}),
\end{align}
\end{subequations}
which defines a $d$-dimensional approximation process $Y^h$ with $Y^h(t_n)=Y_n$.
Here, $\J_{k}$, $k=1,\dots,2m$, are independent random variables which do not depend on $h_n$ and whose moments all exist.
Further, $\alpha=(\alpha_1,\dots,\alpha_s)^\top$, $A=(a_{ij})_{i,j=1,\dots,s}$, $c=(c_1,\dots,c_s)^\top$, $b_{1}=(b_{1,1},\dots,b_{1,s})^\top$, and $b_{2}=(b_{2,1},\dots,b_{2,s})^\top$ are the coefficients of the SRK method. In the following we choose $c=A\Ids$ with $\Ids=(1,\dots,1)^\top\in\mathbb{R}^s$. Consequently, from now on we can assume for the analysis of this methods that SDE \eqref{aco:eq:SDE} is given in autonomous form, i.\,e., $g_0(t,X)\equiv g_0(X)$. The analysis relies on the theory of stochastic B-series, which is shortly reviewed in Section \ref{aco:sec:BSeries} and applied in Section \ref{aco:sec:OrdCond} to derive order conditions for method \eqref{aco:eq:SRK} up to order three. Then, in Section \ref{aco:sec:AN3D1} a concrete explicit third order method is constructed by minimizing the error coefficients. Finally, in Section \ref{aco:sec:numex} we give some numerical examples.
\section{Stochastic B-series}
\label{aco:sec:BSeries}\begin{sloppypar}
Order conditions for method \eqref{aco:eq:SRK} can be calculated using the colored rooted tree theories derived for the weak approximation of It\^{o} respectively Stratonovich SDEs by SRK methods, compare \cite{roessler04ste,roessler06rta,komori07mcr}. Here, we will follow the more general approach developed in \cite{debrabant08bao}, which is based on the work in \cite{burrage96hso,burrage00oco,roessler06rta} and applicable both for It\^{o}- and Stratonovich SDEs as well as strong and weak approximation. For more details and proofs, see \cite{debrabant08bao}.
\end{sloppypar}
First, we introduce the set of colored, rooted trees related to the SDE \eqref{aco:eq:SDE}, as well as the elementary differentials associated with each of these trees. We adapt these definitions to the special case of additive noise by neglecting all terms which are related to derivatives of $g_l$, $l=1,\dots,m$.
\begin{defi}[trees]
    The set of $m+1$-colored, rooted trees \[\Tadd=\{\emptyset\}\cup T_0 \cup \{\bullet_1,\dots,\bullet_m\}\] related to additive noise is recursively defined as follows:
    \begin{enumerate}
    \item[\rm (a)] The graph $\bullet_0=[\emptyset]_0$ with only one vertex
        of color $0$ belongs to $T_0$.
    \end{enumerate}
    Let $\tau=[\tau_1,\tau_2,\ldots,\tau_{\kappa}]_0$ be the tree
    formed by joining the subtrees
    $\tau_1,\tau_2,\ldots,\tau_{\kappa}$ each by a single branch to a
    common root of color $0$.
    \begin{enumerate}
    \item[\rm (b)] If $\tau_1,\tau_2,\ldots,\tau_{\kappa} \in \Tadd$, then
        $\tau=[\tau_1,\tau_2,\ldots,\tau_{\kappa}]_0 \in T_0$.
    \end{enumerate}
\end{defi}
Thus, $T_0$ is the set of trees with a $0$-colored root. $\bullet_0$ will be called deterministic node, $\bullet_l$ for $l>0$ stochastic node of color $l$.
\begin{defi}[elementary differentials]
    For a tree $\tau \in \Tadd$ the elementary differential is
    a mapping $F(\tau):\R^d \rightarrow \R^d$ defined
    recursively by
    \begin{enumerate}
    \item[\rm (a)] $F(\emptyset)(x_0)=x_0$, 
    \item[\rm (b)] $F(\bullet_0)(x_0)=g_0(x_0)$, $F(\bullet_l)(x_0)=g_l$ for $l=1,\dots,m$, 
    \item[\rm (c)] If $\tau=[\tau_1,\tau_2,\ldots,\tau_{\kappa}]_0 \in T_0$,
        then
        \[
        F(\tau)(x_0)=g_0^{(\kappa)}(x_0)
        \left(F(\tau_1)(x_0),F(\tau_2)(x_0),\ldots,F(\tau_{\kappa})(x_0)\right). \]
    \end{enumerate}
\end{defi}
To simplify the presentation, we neglect in the following the index $n$ of $h_n$ and write only $h$. Further, we denote by $\Xi$ the set of families of Borel measurable mappings
\[
\Xi:=\big\{\{\varphi(h)\}_{h\geq0}:\;
    \varphi(h):~\Omega\to\R\text{ is }\mathcal{A}\text{-}\mathcal{B}\text{-measurable }\forall h\geq0
    \big\}.
\]
Both the solution of \eqref{aco:eq:SDE} and its approximation by method \eqref{aco:eq:SRK} can formally be written in terms of B-series.
\begin{defi}[B-series]
    Given a mapping $\phi:\Tadd \rightarrow \Xi$ satisfying
    \[ \phi(\emptyset)\equiv1 \;\text{ and }\; \phi(\tau)(0)=0,\quad
    \forall \tau\in \Tadd \backslash \{\emptyset\}.\]
A (stochastic) B-series is then a formal series of the form
    \[ B(\phi,x_0; h) = \sum_{\tau \in \Tadd}
    \alpha(\tau)\cdot\phi(\tau)(h)\cdot F(\tau)(x_0), \]
where $\alpha: \Tadd\rightarrow \mathbb{Q}$ is given by\vspace{-3pt}
    \begin{align*}
        \alpha(\emptyset)&=1,&\alpha(\bullet_l)&=1,
        &\alpha(\tau=[\tau_1,\ldots,\tau_{\kappa}]_l)&=\frac{1}{r_1!r_2!\cdots
            r_{q}! } \prod_{j=1}^{\kappa} \alpha(\tau_j),
    \end{align*}
    where $r_1,r_2,\ldots,r_{q}$ count equal trees among
    $\tau_1,\tau_2,\ldots,\tau_{\kappa}$.
\end{defi}
For multidimensional $\phi:\Tadd\rightarrow \Xi^{s}$, $s\in\N$, we define \[B(\phi,x_{0};
h)=[B(\phi_{1},x_{0}; h),\ldots, B(\phi_{s},x_{0}; h)]^{\top}.\]
If $Z(h)$ can be written as a B-series,
then $f(Z(h))$ can be written as a similar series, where the sum is
taken over trees with a root of color $f$ and subtrees in $\Tadd$:
\begin{lem} \label{aco:lem:f_y} If $Z(h)=B(\phi, x_0; h)$ is some
    B-series and $f\in C^{\infty}(\R^d,\R^{\hat{d}})$, then
    $f(Z(h))$ can be written as a formal series of the form
    \begin{equation}
        \label{aco:eq:f}
        f(Z(h))=\sum_{u\in U_f^{add}} \beta(u)\cdot \psi_\phi(u)(h)\cdot G(u)(x_0),
    \end{equation}
    where
    \begin{enumerate}
    \item[\rm (a)] $U_f^{add}$ is a set of trees derived from $\Tadd$ as follows: $[\emptyset]_f \in U_f^{add}$, and if
        $\tau_1,\tau_2,\ldots,\tau_{\kappa}\in \Tadd$, then
        $[\tau_1,\tau_2,\ldots,\tau_{\kappa}]_f\in U_f^{add}$,
    \item[\rm (b)] $G([\emptyset]_f)(x_0)=f(x_0)$ and $G(u=[\tau_1,\ldots,\tau_{\kappa}]_f)(x_0) =
        f^{(\kappa)}(x_0)
        \big(F(\tau_1)(x_0),\ldots,F(\tau_{\kappa})(x_0)\big)$,
    \item[\rm (c)] $\beta([\emptyset]_f)=1$ and $
        \beta(u=[\tau_1,\ldots,\tau_{\kappa}]_f)
        =\frac{1}{r_1!r_2!\ldots r_{q}!}\prod_{j=1}^{\kappa}
        \alpha(\tau_{{j}})$, where $r_1,r_2,\ldots,\break r_{q}$ count
        equal trees among $\tau_1,\tau_2,\ldots,\tau_{\kappa}$, 
    \item[\rm (d)] $\psi_\phi([\emptyset]_f)\equiv1$ and
        $\psi_\phi(u=[\tau_1,\ldots,\tau_{\kappa}]_f)(h) =
        \prod_{j=1}^{\kappa} \phi(\tau_j)(h)$.
    \end{enumerate}
\end{lem}
\begin{remark}\rm
To simplify the presentation, we assume throughout this article that all derivatives of $f$ and $g_0$ exist. Otherwise, one had to consider truncated B-series with a remainder term.
\end{remark}
\begin{theo} \label{aco:thm:Bex} The solution $X(t_0+h)$ of \eqref{aco:eq:SDE}
    can be written as a B-series $B(\varphi,x_0; h)$ with
    \begin{align*}
        \varphi(\emptyset)\equiv1,\qquad\varphi(\bullet_0)(h)=h,\qquad\varphi(\bullet_l)(h)=W_l(h),\quad l=1,\dots,m,\\
        \varphi(\tau=[\tau_1,\ldots,\tau_{\kappa}]_0)(h)=\int_0^h
        \prod_{j=1}^{\kappa} \varphi(\tau_j)(s)~ ds.
    \end{align*}
\end{theo}
The following definition of the
order of the tree, $\rho(\tau)$, is motivated by the fact that $\E
W_{l}(h)^{2}=h$ for $l\geq 1 $.
\begin{defi}[order]
    The order of a tree $\tau \in \Tadd$ is defined by
    \[\rho(\emptyset)=0,
    \quad\rho(\bullet_l)=\frac12,~l=1,\dots,m,\]
    and
    \[\rho(\tau=[\tau_1,\dots,\tau_\kappa]_l)=\sum\limits_{i=1}^\kappa\rho(\tau_i)+
    \begin{cases}
        1&\text{for }l=0,\\
        \frac12&\text{for }l>0.
    \end{cases}
    \]
    The order of a tree $u \in U_f^{add}$ is given by $\rho(u=[\tau_1,\dots,\tau_\kappa]_f)=\sum\limits_{i=1}^\kappa\rho(\tau_i)$.
\end{defi}

In the following we define the product of vectors by componentwise multiplication.

\begin{theo} \label{aco:thm:Bnum}The numerical approximation $Y_1$ as well as the stage values can be written in terms of B-series
    \[ H = B\left(\PhiH,x_0; h\right), \qquad Y_1 =
    B(\Phi,x_0; h)\]
with
    \begin{subequations} \label{aco:eq:B_iter_H}
        \begin{eqnarray}
            &\PhiH(\emptyset)\equiv \Ids, \quad
            \PhiH(\bullet_l)(h)=\sqrt{h}(b_1\J_l+b_2\J_{m+l}),\quad l=1,\dots,m,& \\
            &\PhiH(\tau=[\tau_1,\ldots,\tau_{\kappa}]_0)(h)
            = h\A\prod_{j=1}^{\kappa}
            \PhiH(\tau_j)(h)&
        \end{eqnarray}
    \end{subequations}
    and
    \begin{subequations} \label{aco:eq:B_iter_Y}
        \begin{eqnarray}
            &\Phi(\emptyset)\equiv 1, \quad
            \Phi(\bullet_l)(h)=\sqrt{h}\J_l,\quad l=1,\dots,m,& \\
            &\Phi(\tau=[\tau_1,\ldots,\tau_{\kappa}]_0)(h)
            = h\alpha^\top\prod_{j=1}^{\kappa}
            \PhiH(\tau_j)(h).&
        \end{eqnarray}
    \end{subequations}
\end{theo}
\section{Derivation of order conditions}
\label{aco:sec:OrdCond}
With all the B-series in place, we can now present the order
conditions for the weak convergence.

Let $le_f(h;t,x)$ be the weak local error of the method starting
at the point $(t,x)$ with respect to the functional $f$ and step size
$h$, i.\,e.
\[
le_f(h;t,x)= \E \big(f(Y^h(t+h))-f(X(t+h))|Y^h(t)=X(t)=x\big).
\]
From Theorems \ref{aco:thm:Bex} and \ref{aco:thm:Bnum} and Lemma
\ref{aco:lem:f_y} we obtain
\begin{equation}\label{aco:eq:lokError}
le_f(h;t,x)=\sum_{u\in
    U_f^{add}}\beta(u)\cdot\E\left[\psi_\Phi(u)(h)-\psi_\varphi(u)(h)\right]\cdot
  G(u)(x)
\end{equation}
with
\begin{equation}\label{aco:eq:psivarphi}
\psi_\varphi([\emptyset]_f)\equiv1,\quad
\psi_\varphi(u=[\tau_1,\ldots,\tau_\kappa]_f)(h) =
\prod\limits_{j=1}^\kappa \varphi(\tau_j)(h)
\end{equation}
and
\begin{equation}\label{aco:eq:psiPhi}
\psi_\Phi([\emptyset]_f)\equiv1,\quad
\psi_\Phi(u=[\tau_1,\ldots,\tau_\kappa]_f)(h) =
\prod\limits_{j=1}^\kappa \Phi(\tau_j)(h).
\end{equation} Thus, we have weak consistency of order $p$ (and thus, due to the Milstein theorem \cite{milstein95nio}, also weak convergence) if and only if
\begin{equation}
    \label{aco:eq:T_eq}
    \E\psi_\Phi(u)(h)=\E\psi_\varphi(u)(h)+\bO(h^{p+1})\quad\forall u\in U_f^{add}\text{ with }\rho(u)\leq p+\frac12.
\end{equation}
Note that \eqref{aco:eq:T_eq} slightly weakens conditions given in \cite{roessler06rta}.

By Theorems \ref{aco:thm:Bex} and \ref{aco:thm:Bnum}, \eqref{aco:eq:psivarphi} and \eqref{aco:eq:psiPhi} we can now evaluate the order conditions \eqref{aco:eq:T_eq} and obtain the following theorem.

\begin{theo}\label{aco:thm:OrderCond}For a $p$-th order method choose the independent random variables $\J_k$ of the SRK method \eqref{aco:eq:SRK} such that their moments coincide with those of $N(0,1)$ up to the $(2p+1)$-th moment for k=1,\dots,m and up to the $(2p-1)$-th moment for k=m+1,\dots,2m.
If in addition the coefficients of the SRK method \eqref{aco:eq:SRK} fulfill
\begin{multienumerate}
\mitemx{$\alpha^\top\Ids=1$,}
\end{multienumerate}
then the method is of weak order $p=1$. If also the equations
\begin{multienumerate}
\addtocounter{multienumi}{1}
\mitemxxx{$\alpha^\top A\Ids=\frac12,$}{$\alpha^\top(b_1^2+b_2^2)=\frac12,$}{$\alpha^\top b_1=\frac12$}
\end{multienumerate}
are fulfilled, then the SRK method is of weak order $p=2$. Finally, if additionally
\begin{multienumerate}
\addtocounter{multienumi}{4}
\mitemxxx{$\alpha^\top A^2\Ids=\frac16,$}
{$\alpha^\top(A\Ids)^2=\frac13,$}
{$\alpha^\top A(b_1^2+b_2^2)=\frac16,$}
\mitemxox{$\alpha^\top(b_1(Ab_1)+b_2(Ab_2))=\frac16,$}
{$\alpha^\top Ab_1=\frac16,$}
\mitemxox
{$\alpha^\top((A\Ids)(b_1^2+b_2^2))=\frac13,$}
{$\alpha^\top((A\Ids)b_1)=\frac13,$}
\mitemxxx{$\alpha^\top(b_1^2+b_2^2)^2=\frac13,$}
{$\alpha^\top(b_1^3+b_1b_2^2)=\frac13,$}
{$\alpha^\top b_1^2=\frac13,$}
\mitemx{$(\alpha^\top b_2)^2=\frac1{12}$}
\end{multienumerate}
are fulfilled, then the SRK method is of weak order $p=3$.
\end{theo}
\begin{proof}
First, we note that $\E\psi_\varphi(u)=0$ for all trees $u\in U_f^{add}$ which have an odd number of stochastic nodes of one color, see \cite{debrabant10ste} or also \cite{burrage99rkm,roessler04ste}.
For those of these trees which have an order $\rho(u)\leq p+\frac12$, by construction of the method and due to the assumptions on $\J_k$, $k=1,\dots,2m$, it holds also
$\E\psi_\Phi(u)=0$. Thus, in the following we only have to consider trees with an even number of each kind of stochastic nodes, in particular only trees of integer order.
Consequently, there are only two kinds of trees of order one to consider:
\[
u_1=\raisebox{-0.3cm}{\scalebox{0.7}{
\pstree{\Tcircle{f}}{
\Tcircle{j}\Tcircle{j}}}}
,\quad j=1,\dots,m,
\qquad\text{and}\qquad
u_2=\raisebox{-0.3cm}{\scalebox{0.7}{
\pstree{\Tcircle{f}}{
\Tcircle{0}}}}.
\]
Theorems \ref{aco:thm:Bex} and \ref{aco:thm:Bnum}, \eqref{aco:eq:psivarphi} and \eqref{aco:eq:psiPhi} yield
\[\psi_\Phi(u_1)=h\J_j^2,\quad\psi_\varphi(u_1)=W_j(h)^2,\quad \psi_\Phi(u_2)=h\alpha^\top\Ids\quad\text{and}\quad\psi_\varphi(u_2)=h.\] Thus, by the assumptions on $\J_j$, $\E\psi_\Phi(u_1)(h)=\E\psi_\varphi(u_1)(h)$ is fulfilled automatically, whereas $\E\psi_\Phi(u_2)(h)=\E\psi_\varphi(u_2)(h)$ yields order condition 1.

If $u\in U_f^{add}$ with $u=[\tau_1,\dots,\tau_\kappa]_f$ can be split into two trees $u_1=[\tau_{i_1},\dots,\tau_{i_{\kappa_1}}]_f$,
$u_2=[\tau_{j_1},\dots,\tau_{j_{\kappa_2}}]_f$ with disjoint stochastic nodes, i.\,e.\ such that $\kappa_1,\kappa_2>0$, $\kappa_1+\kappa_2=\kappa$, $\{i_1,\dots,i_{\kappa_1},j_1,\dots,j_{\kappa_2}\}=\{1,\dots,\kappa\}$, and the sets of colors of the stochastic nodes of $u_1$ and $u_2$ are disjoint, then
\[
\E\psi_\Phi(u)(h)=\E\psi_\Phi(u_1)(h)\E\psi_\Phi(u_2)(h)
=\E\psi_\varphi(u_1)(h)\E\psi_\varphi(u_2)(h)=\E\psi_\varphi(u)(h),
\]
provided that the order conditions of orders lower than $\rho(u)$ are fulfilled. Thus, in the following we only have to consider trees of second and third order which cannot be decomposed into two trees with disjoint stochastic nodes. The relevant second order trees together with the derivation of the corresponding order conditions are given in Table \ref{aco:tab:reltreesord2}, the ones of order three in Tables \ref{aco:tab:reltreesord3a}-\ref{aco:tab:reltreesord3c}, which completes the proof.
\newcommand{\bz}[6]{$\begin{array}{c}
\\[-2mm]#1\\[-2mm]\\
\end{array}$&$\begin{array}{c}
#2\\#3
\end{array}$&$\begin{array}{c}
#4\\#5
\end{array}$&#6\\}
\newcommand{\bzt}[7]{\bz{\raisebox{#7}{\scalebox{0.7}{#1}}}{#2}{#3}{#4}{#5}{#6}}
\newcommand{\bztd}[6]{\bzt{#1}{#2}{#3}{#4}{#5}{#6}{-0.6cm}}
\newcommand{\bztz}[6]{\bzt{#1}{#2}{#3}{#4}{#5}{#6}{-0.3cm}}
\begin{table}
\caption{\label{aco:tab:reltreesord2}Relevant second order trees and derivation of corresponding order conditions}
\begin{tabular}{c|c|c|c}
\hline\noalign{\smallskip}
\bz{u}{\psi_\Phi(u)}{\E\psi_\Phi(u)}{\psi_\varphi(u)}{\E\psi_\varphi(u)}{ord. cond.}
\noalign{\smallskip}\hline\noalign{\smallskip}
\bztz{\pstree{\Tcircle{f}}{\Tcircle{j}\Tcircle{j}\Tcircle{j}\Tcircle{j}}}
{h^2\J_j^4}{3h^2}{W_j(h)^4}{3h^2}{by assumption}
\hline
\bztd{\pstree{\Tcircle{f}}{\pstree{\Tcircle{0}}{
\Tcircle{0}}}}{h^2\alpha^\top A\Ids}{h^2\alpha^\top A\Ids}{\int_0^hs~ds}{\frac{h^2}2}{2.}
\hline
\bztd{\pstree{\Tcircle{f}}{\pstree{\Tcircle{0}}{
\Tcircle{j}\Tcircle{j}}}}{h^2\alpha^\top(b_1\J_j+b_2\J_{m+j})^2}{h^2\alpha^\top(b_1^2+b_2^2)}
{\int_0^hW_j(s)^2~ds}{\frac{h^2}2}{3.}
\hline
\bztd{\pstree{\Tcircle{f}}{\Tcircle{j}\pstree{\Tcircle{0}}{
\Tcircle{j}}}}{h^2\J_j\alpha^\top(b_1\J_j+b_2\J_{m+j})}{h^2\alpha^\top b_1}{W_j(h)\int_0^hW_j(s)~ds}{\frac{h^2}2}{4.}
\noalign{\smallskip}\hline
\end{tabular}
\end{table}
\begin{subtables}
\newcommand{\bztv}[6]{\bzt{#1}{#2}{#3}{#4}{#5}{#6}{-0.9cm}}
\begin{table}
\caption{\label{aco:tab:reltreesord3a}Relevant third order trees and derivation of corresponding order conditions, part 1}
\begin{tabular}{c|c|c|c@{}}
\hline\noalign{\smallskip}
\bz{u}{\psi_\Phi(u)}{\E\psi_\Phi(u)}{\psi_\varphi(u)}{\E\psi_\varphi(u)}{ord. cond.}\noalign{\smallskip}\hline\noalign{\smallskip}
\bztz{\pstree{\Tcircle{f}}{\Tcircle{j}\Tcircle{j}\Tcircle{j}\Tcircle{j}\Tcircle{j}\Tcircle{j}}}{h^3\J_j^6}{15h^2}{W_j(h)^6}{15h^2}{by assumption}\hline
\bztv{\pstree{\Tcircle{f}}{\pstree{\Tcircle{0}}{\pstree{\Tcircle{0}}{
\Tcircle{0}}}}}{h^3\alpha^\top A^2\Ids}{h^3\alpha^\top A^2\Ids}{\int_0^h\int_0^{s_1}s_2~ds_2~ds_1}{\frac{h^3}6}{5.}
\hline
\bztd{\pstree{\Tcircle{f}}{\pstree{\Tcircle{0}}{
\Tcircle{0}\Tcircle{0}}}}{h^3\alpha^\top(A\Ids)^2}{h^3\alpha^\top (A\Ids)^2}{\int_0^hs^2~ds}{\frac{h^3}3}{6.}
\hline
\bztv{\pstree{\Tcircle{f}}{\pstree{\Tcircle{0}}
{\pstree{\Tcircle{0}}{\Tcircle{j}\Tcircle{j}}}}}{h^3\alpha^\top A(b_1\J_j+b_2\J_{m+j})^2}{h^3\alpha^\top A(b_1^2+b_2^2)}{\int_0^h\int_0^{s_1}W_j(s_2)^2~ds_2~ds_1}{\frac{h^3}6}{7.}
\noalign{\smallskip}\hline
\end{tabular}
\end{table}
\begin{table}
\caption{\label{aco:tab:reltreesord3b}Relevant third order trees and derivation of corresponding order conditions, part 2}
\begin{tabular}{c|@{}c@{}|@{}c@{}@{}|@{}c@{}}
\hline\noalign{\smallskip}
\bz{u}{\psi_\Phi(u)}{\E\psi_\Phi(u)}{\psi_\varphi(u)}{\E\psi_\varphi(u)}{\begin{tabular}{c}
ord.\\
cond.
\end{tabular}}
\noalign{\smallskip}\hline\noalign{\smallskip}
\bztv{\pstree{\Tcircle{f}}{\pstree{\Tcircle{0}}
{\Tcircle{j}\pstree{\Tcircle{0}}{\Tcircle{j}}}}}{h^3\alpha^\top ((b_1\J_j+b_2\J_{m+j})(A(b_1\J_j+b_2\J_{m+j})))}{h^3\alpha^\top (b_1(Ab_1)+b_2(Ab_2))}{\int_0^hW_j(s_1)\int_0^{s_1}W_j(s_2)~ds_2~ds_1}{\frac{h^3}6}{8.}
\hline
\bztv{\pstree{\Tcircle{f}}{\Tcircle{j}\pstree{\Tcircle{0}}
{\pstree{\Tcircle{0}}{\Tcircle{j}}}}}{h^3\J_j\alpha^\top A(b_1\J_j+b_2\J_{m+j})}{h^3\alpha^\top Ab_1}{W_j(h)\int_0^h\int_0^{s_1}W_j(s_2)~ds_2~ds_1}{\frac{h^3}6}{9.}
\hline
\bztd{\pstree{\Tcircle{f}}{\pstree{\Tcircle{0}}{
\Tcircle{0}\Tcircle{j}\Tcircle{j}}}}{h^3\alpha^\top((A\Ids)(b_1\J_j+b_2\J_{m+j})^2)}
{h^3\alpha^\top((A\Ids)(b_1^2+b_2^2))}{\int_0^hsW_j(s)^2~ds}{\frac{h^3}3}{10.}
\hline
\bztd{\pstree{\Tcircle{f}}{\Tcircle{j}\pstree{\Tcircle{0}}{
\Tcircle{0}\Tcircle{j}}}}{h^3\J_j\alpha^\top((A\Ids)(b_1\J_j+b_2\J_{m+j}))}
{h^3\alpha^\top((A\Ids)b_1)}{W_j(h)\int_0^hsW_j(s)~ds}{\frac{h^3}3}{11.}
\hline
\bztd{\pstree{\Tcircle{f}}{\pstree{\Tcircle{0}}{
\Tcircle{j}\Tcircle{j}\Tcircle{k}\Tcircle{k}}}}{h^3\alpha^\top((b_1\J_j+b_2\J_{m+j})^2(b_1\J_k+b_2\J_{m+k})^2)}{\begin{cases}
h^3\alpha^\top(3b_1^4+6b_1^2b_2^2+3b_2^4)&j=k\\
h^3\alpha^\top(b_1^2+b_2^2)^2&j\neq k
\end{cases}
}
{\int_0^hW_j(s)^2W_k(s)^2~ds}{\begin{cases}
h^3&j=k\\
\frac{h^3}3&j\neq k
\end{cases}}{12.}
\noalign{\smallskip}\hline
\end{tabular}
\end{table}
\begin{table}
\caption{\label{aco:tab:reltreesord3c}Relevant third order trees and derivation of corresponding order conditions, part 3}
\begin{tabular}{@{}c|@{}c@{}|@{}c@{}|@{}c@{}}
\hline\noalign{\smallskip}
\bz{u}{\psi_\Phi(u)}{\E\psi_\Phi(u)}{\psi_\varphi(u)}{\E\psi_\varphi(u)}{
\begin{tabular}{c}
ord.\\
cond.
\end{tabular}}
\noalign{\smallskip}\hline\noalign{\smallskip}
\bztd{\pstree{\Tcircle{f}}{\Tcircle{k}\pstree{\Tcircle{0}}{
\Tcircle{j}\Tcircle{j}\Tcircle{k}}}}{h^3\J_k\alpha^\top((b_1\J_j+b_2\J_{m+j})^2(b_1\J_k+b_2\J_{m+k}))}{\begin{cases}
h^3\alpha^\top(3b_1^3+3b_1b_2^2)&j=k\\
h^3\alpha^\top(b_1^3+b_1b_2^2)&j\neq k
\end{cases}
}
{W_k(h)\int_0^hW_j(s)^2W_k(s)~ds}{\begin{cases}
h^3&j=k\\
\frac{h^3}3&j\neq k
\end{cases}}{13.}
\hline
\bztd{\pstree{\Tcircle{f}}{\Tcircle{j}\Tcircle{k}\pstree{\Tcircle{0}}{
\Tcircle{j}\Tcircle{k}}}}{h^3\J_j\J_k\alpha^\top((b_1\J_j+b_2\J_{m+j})(b_1\J_k+b_2\J_{m+k}))}
{\begin{cases}
h^3\alpha^\top(3b_1^2+b_2^2)&j=k\\
h^3\alpha^\top b_1^2&j\neq k
\end{cases}}{W_j(h)W_k(h)\int_0^hW_j(s)W_k(s)~ds}
{\begin{cases}\frac{7h^3}6&j=k\\
\frac{h^3}3&j\neq k
\end{cases}}{\begin{tabular}{c}
due\\ to 3.\\
equiv.\\ to 14.
\end{tabular}}
\hline
\bztd{\pstree{\Tcircle{f}}{\Tcircle{j}\Tcircle{j}\Tcircle{j}\pstree{\Tcircle{0}}{
\Tcircle{j}}}}{h^3\J_j^3\alpha^\top(b_1\J_j+b_2\J_{m+j})}
{3h^3\alpha^\top b_1}{W_j(h)^3\int_0^hW_j(s)~ds}{\frac{3h^3}2}{4.}
\hline
\bztd{\pstree{\Tcircle{f}}{\pstree{\Tcircle{0}}{\Tcircle{j}}\pstree{\Tcircle{0}}{\Tcircle{j}}}}{h^3(\alpha^\top(b_1\J_j+b_2\J_{m+j}))^2}
{h^3((\alpha^\top b_1)^2+(\alpha^\top b_2)^2)}{(\int_0^hW_j(s)~ds)^2}{\frac{h^3}3}{\begin{tabular}{c}
due\\ to 4.\\
equiv.\\ to 15.
\end{tabular}}
\noalign{\smallskip}\hline
\end{tabular}
\end{table}\label{aco:tab:reltreesord3}
\end{subtables}
\end{proof}
Possible discrete choices for the random variables $J_k$, $k=1,\dots,2m$, can be found in Table \ref{aco:Table:DiskrRandomVariables}.
\begin{table}
\caption{\label{aco:Table:DiskrRandomVariables} Some discrete random variables corresponding up to the $i$th moment to $N(0,1)$}
$\begin{array}{c|l}
\hline\noalign{\smallskip}
i&distribution\\\noalign{\smallskip}\hline\noalign{\smallskip}
1&P(\J_k=0)=1\\
3&P(\J_k=1)=P(\J_k=-1)=\frac12\\
5&P(\J_k=\sqrt3)=P(\J_k=-\sqrt3)=\frac16,\quad P(\J_k=0)=\frac23\\
7&P(\J_k=\sqrt6)=P(\J_k=-\sqrt6)=\frac1{30},\quad
P(\J_k=1)=P(\J_k=-1)=\frac3{10},\quad
P(\J_k=0)=\frac13
\\\noalign{\smallskip}\hline
\end{array}$
\end{table}
\section{A concrete explicit third order SRK method}
\label{aco:sec:AN3D1}
Based on Theorem \ref{aco:thm:OrderCond}, we now calculate the coefficients of
an explicit third order SRK method.
The coefficients will be arranged in
an extended Butcher array of the form
\renewcommand{\arraystretch}{1.8}
\[
\begin{array}{c|c|c|c}
    c & A&b_1&b_2\\
    \cline{1-4}
    \hline
    & \alpha^T
\end{array}.
\]
\renewcommand{\arraystretch}{1.0}
Whereas in the deterministic case we would only need three stages to construct an explicit third order method, here we need four stages to fulfill the 15 order conditions of Theorem \ref{aco:thm:OrderCond}. Therefore, we consider $s=4$ in \eqref{aco:eq:SRK}, but require in addition that the method fulfills also the deterministic order four conditions. The remaining degrees of freedom are then eliminated by minimizing the vector $lec$ of the order four coefficients of the local error \eqref{aco:eq:lokError} in the Euclidean norm assuming two dimensional noise ($m=2$), i.\,e.\ by minimizing $\|lec\|_2$ where
\[
lec=\left(\beta(u)\cdot\E\left[\psi_\Phi(u)(h)-\psi_\varphi(u)(h)\right]\right)_{
u\in U_f^{add},
\rho(u)=4
}
.\]
Using again the B-series analysis, a tedious calculation
(one obtains 52 non automatically vanishing terms)
and a subsequent attempt of numerical optimization yield the scheme
AN3D1 presented in Table \ref{aco:Table-Coeff-AN3D1}.
\begin{table}
\caption{\label{aco:Table-Coeff-AN3D1}Coefficients of AN3D1}
$\begin{array}{c|cccc|c|c}
0&0&0&0&0&\multirow{4}{*}{$b_1$}&\multirow{4}{*}{$b_2$}\\
1&1&0&0&0&\\
1/2&3/8&1/8&0&0&\\
1&-0.4526683126055039&-0.4842227708685013&1.9368910834740051&0&\\
    \cline{1-4}
    \hline
    &1/6&-0.005430430675258792&2/3&0.1720970973419255
\end{array}$
with
\begin{align*}
b_1&=(-0.01844540496323970,0.8017012756521233,0.5092227024816198,0.9758794209767762)^\top\\
b_2&=(-0.1866426386543421,-0.8575745885712401,-0.4723392695015512,0.3060354860326548)^\top
\end{align*}
\end{table}
AN3D1 needs two random variable and four drift evaluations per step, and thus two random variable and three drift evaluations less than Platen's third order method.
\section{Numerical example}
\label{aco:sec:numex}
In the following we compare for three simple test equations the performance of the SRK scheme AN3D1 (with N(0,1)-distributed random variables) presented in the last section with some well known schemes, namely the third and the second order SRK schemes due to Platen~\cite{kloeden99nso}, denoted here by PL3 and PL2, respectively, DRI1 due to Debrabant and R\"{o}{\ss}ler~\cite{debrabant09foe}, and the ex\-tra\-po\-la\-ted Euler-Maruyama scheme EXEM (cp.\ \cite{talay90eot}) also attaining order two, which is given by $2\E(\left(Z^{h/2}(t)\right)^2)-\E(\left(Z^{h}(t)\right)^2)$, based on the Euler-Maruyama approximations $Z^{h/2}(t)$ and $Z^h(t)$ calculated with step sizes $h/2$ and $h$.
In each case, the functional $u(t)=\E(f(X(t)))$ is approximated by a Monte Carlo simulation. The sample average $u_{M,h}(t) = \frac{1}{M} \sum_{k=1}^Mf\left(Y^h(t,\omega_k)\right)$, $\omega_k \in \Omega$, of $M=10^9$ independent simulated realizations of the considered approximation $Y^h(t)$ is calculated in order to estimate the expectation and thus to determine the systematic error of the considered schemes. In the following, we denote by $\hat{\mu} = u_{M,h}(T) - u(T)$ the mean error at time $T$ and by $\hat{\sigma}^2_{\mu}$ the empirical variance of the mean error.
Further, we calculate the confidence interval with boundaries $a$
and $b$ to the level of 90\% for the estimated error $\hat{\mu}$
(see \cite{kloeden99nso} for details).

\begin{figure}[tbp]
\includegraphics[width=6.8cm]{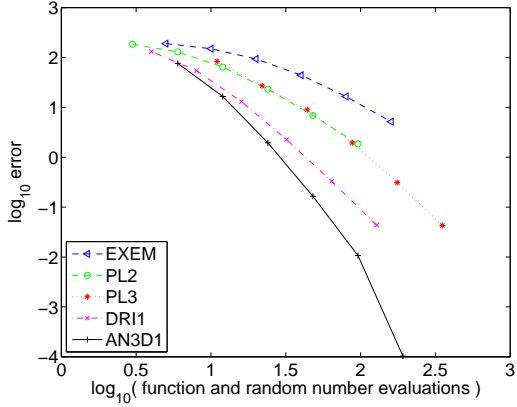}
\caption{Computational effort per simulation path versus precision for SDE \eqref{eq:ex1}} \label{aco:Bild001}
\end{figure}

First, we compute the second moment of the solution of the linear SDE
\begin{equation}\label{eq:ex1}
X(t)=\frac1{10}+\frac32\int_0^tX(s)~ds+\frac1{10}W(t),\quad t\in I=[0,2],
\end{equation}
which can be calculated analytically as
\begin{equation}
    \E(X^2(t)) = \frac29\left(\frac{397}{200}-\frac{23}5e^{3/2t}+\frac{133}{50}e^{3t}\right).
\end{equation}
The solution value $\E(X^2(T))$ is now approximated with step sizes $2^{1},\ldots, 2^{-4}$  at time $T=2$.
The results for the applied schemes are presented in Table~\ref{aco:Table1}.
Of course, these results have to be related to the computational effort of the schemes which we take in the following as sum of the number of evaluations of the drift function $a$  as well as the number of random variables that have to be simulated. Then we can oppose the computational efforts to the errors of the analyzed schemes.
The results are presented in Figure~\ref{aco:Bild001}. Although being of different order, the two Platen schemes yield comparable results. This is due to the much higher computational costs of PL3. Both methods are better than the extrapolated Euler method, but perform worse than DRI1, which has optimized coefficients \cite{debrabant09foe} and behaves therefore nearly like an order three method. Our new method AN3D1 performs best.
\begin{table}
\caption{Mean errors, empirical variances and confidence intervals for SDE \eqref{eq:ex1}} \label{aco:Table1}
\begin{center}
\begin{tabular}{c|c|c|c|c|c}
   & $h$ & $\hat{\mu}$ & $\hat{\sigma}_{\mu}^2$  & $a$ & $b$ \\
   \hline
\multirow{4}{*}{EXEM}
& $2^{1}$ &-1.900E+02 & 7.882E-07 & -1.900E+02 & -1.900E+02\\
& $2^{0}$ &-1.499E+02 & 7.032E-06 & -1.499E+02 & -1.499E+02\\
& $2^{-1}$ &-9.357E+01 & 3.666E-05 & -9.357E+01 & -9.357E+01\\
& $2^{-2}$ &-4.435E+01 & 8.881E-05 & -4.435E+01 & -4.434E+01\\
& $2^{-3}$ &-1.649E+01 & 1.988E-04 & -1.650E+01 & -1.649E+01\\
& $2^{-4}$ & -5.170E+00 & 2.770E-04 & -5.174E+00 & -5.166E+00\\
\hline
\multirow{4}{*}{PL2}
& $2^{1}$ &-1.840E+02 & 6.783E-07 & -1.840E+02 & -1.840E+02\\
& $2^{0}$ &-1.294E+02 & 6.348E-06 & -1.294E+02 & -1.294E+02\\
& $2^{-1}$ &-6.412E+01 & 2.790E-05 & -6.412E+01 & -6.412E+01\\
& $2^{-2}$ &-2.312E+01 & 4.995E-05 & -2.312E+01 & -2.312E+01
\\
& $2^{-3}$ &-6.880E+00 & 5.861E-05 & -6.882E+00 & -6.878E+00\\
& $2^{-4}$ &-1.863E+00 & 8.264E-05 & -1.865E+00 & -1.861E+00\\
       \hline
\multirow{4}{*}{PL3}
& $2^{1}$ &-8.377E+01 & 2.936E-03 & -8.378E+01 & -8.375E+01\\
& $2^{0}$ &-2.705E+01 & 1.614E-03 & -2.706E+01 & -2.704E+01\\
& $2^{-1}$ &-8.941E+00 & 2.345E-04 & -8.944E+00 & -8.937E+00\\
& $2^{-2}$ &-1.951E+00 & 6.624E-05 & -1.953E+00 & -1.949E+00\\
& $2^{-3}$ &-3.111E-01 & 6.180E-05 & -3.130E-01 & -3.093E-01\\
& $2^{-4}$ &-4.307E-02 & 4.718E-05 & -4.470E-02 & -4.144E-02\\
\hline
\multirow{4}{*}{DRI1}
& $2^{1}$ &-1.316E+02 & 3.486E-06 & -1.316E+02 & -1.316E+02\\
& $2^{0}$ &-5.438E+01 & 2.049E-05 & -5.438E+01 & -5.437E+01\\
& $2^{-1}$ &-1.308E+01 & 4.872E-05 & -1.308E+01 & -1.308E+01\\
& $2^{-2}$ &-2.254E+00 & 6.097E-05 & -2.256E+00 & -2.252E+00 \\
& $2^{-3}$&-3.314E-01 & 6.225E-05 & -3.333E-01 & -3.295E-01\\
& $2^{-4}$&-4.343E-02 & 8.405E-05 & -4.560E-02 & -4.126E-02\\
\hline
\multirow{4}{*}{AN3D1}
& $2^{1}$ &-7.638E+01 & 1.286E-05 & -7.638E+01 & -7.638E+01\\
& $2^{0}$ &-1.654E+01 & 4.729E-05 & -1.654E+01 & -1.654E+01\\
& $2^{-1}$ &-1.946E+00 & 6.804E-05 & -1.948E+00 & -1.944E+00\\
& $2^{-2}$ &-1.651E-01 & 4.993E-05 & -1.668E-01 & -1.635E-01\\
& $2^{-3}$ &-1.073E-02 & 5.940E-05 & -1.255E-02 & -8.900E-03\\
& $2^{-4}$&-1.030E-04 & 4.754E-05 & -1.738E-03 & 1.532E-03
\end{tabular}
\end{center}
\end{table}

\begin{figure}[tbp]
\includegraphics[width=6.8cm]{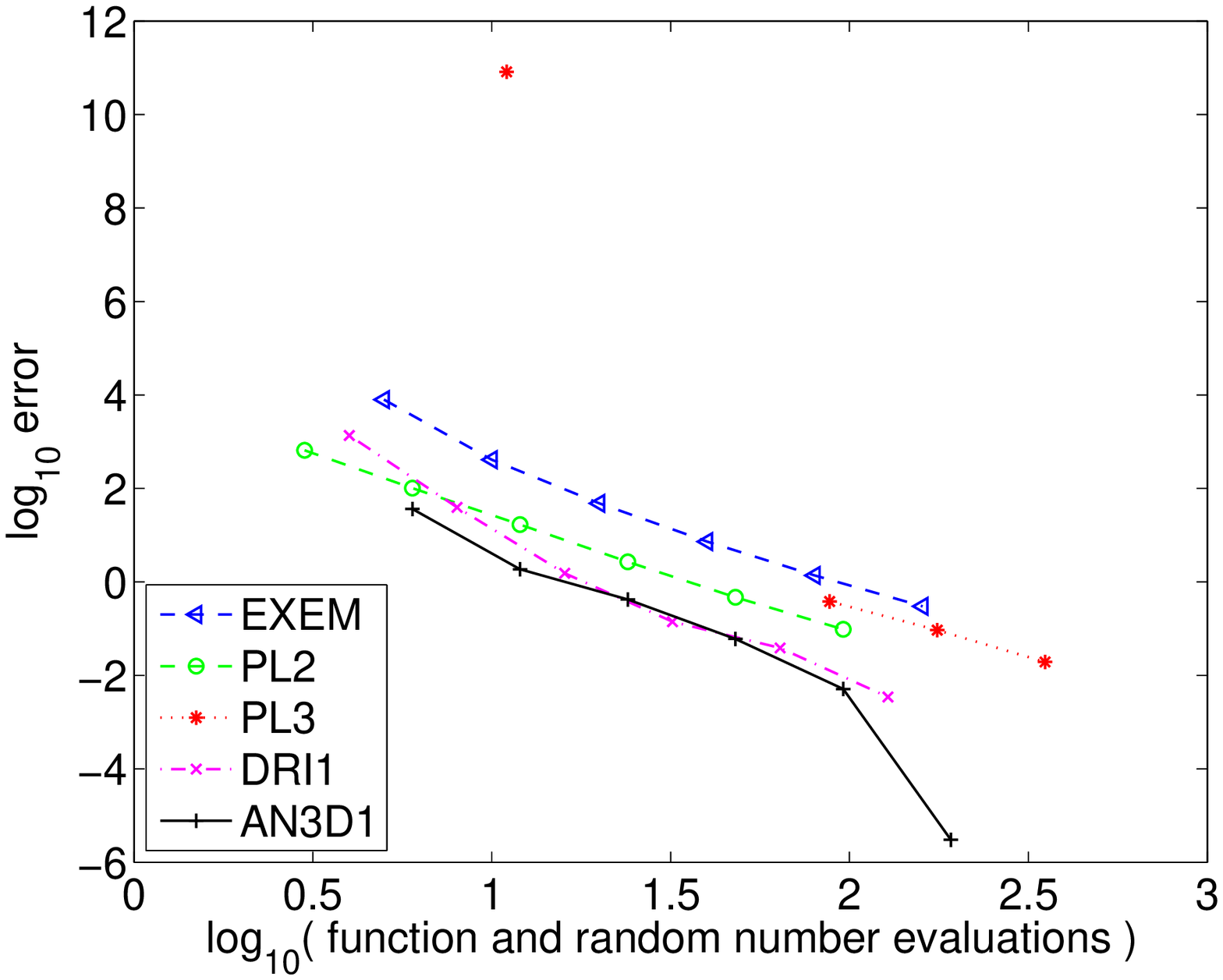}
\caption{Computational effort per simulation path versus precision for SDE \eqref{eq:ex2}} \label{aco:Bild002}
\end{figure}

As next example we consider the nonlinear SDE
\begin{equation}\label{eq:ex2}
X(t)=\frac1{10}+\int_0^t\left(\frac32e^{-2X(s)}+1\right)~ds+\frac1{10} W(t),\quad t\in I=[0,2].
\end{equation}
Then $\E(e^{2X(t)})$ can be calculated as
\begin{equation}
    \E(e^{2X(t)})) =\left(e^{1/5}+\frac{150}{101}\right)e^{101/50 t}-\frac{150}{101}.
\end{equation}
The solution value $\E(e^{2X(T)})$ is approximated with step sizes $2^{1},\ldots, 2^{-4}$ at time $T=2$.
The results for the applied schemes are presented in Table~\ref{aco:Table2} and Figure~\ref{aco:Bild002} and reflect a similar behaviour to the one from the linear example, except that PL3 suffers now from stability problems.
\begin{table}
\caption{Mean errors, empirical variances and confidence intervals for SDE \eqref{eq:ex2}} \label{aco:Table2}
\begin{center}
\begin{tabular}{c|c|c|c|c|c}
   & $h$ & $\hat{\mu}$ & $\hat{\sigma}_{\mu}^2$  & $a$ & $b$ \\
   \hline
\multirow{4}{*}{EXEM}
& $2^{1} $ &-7.925E+03 & 2.818E-01 & -7.925E+03 & -7.925E+03\\
& $2^{0} $ &-4.127E+02 & 2.748E-03 & -4.127E+02 & -4.127E+02\\
& $2^{-1}$ &-4.777E+01 & 8.555E-04 & -4.777E+01 & -4.776E+01\\
& $2^{-2}$ &-7.296E+00 & 6.418E-04 & -7.302E+00 & -7.290E+00\\
& $2^{-3}$ &-1.369E+00 & 5.752E-04 & -1.375E+00 & -1.363E+00\\
& $2^{-4}$ &-2.992E-01 & 5.557E-04 & -3.048E-01 & -2.936E-01\\
\hline
\multirow{4}{*}{PL2}
& $2^{1} $&6.573E+02 & 2.146E-03 & 6.573E+02 & 6.573E+02\\
& $2^{0} $&1.010E+02 & 1.921E-04 & 1.010E+02 & 1.011E+02\\
& $2^{-1}$&1.678E+01 & 1.192E-04 & 1.678E+01 & 1.678E+01\\
& $2^{-2}$&2.676E+00 & 1.064E-04 & 2.674E+00 & 2.679E+00\\
& $2^{-3}$&4.665E-01 & 1.064E-04 & 4.641E-01 & 4.690E-01\\
& $2^{-4}$&9.702E-02 & 1.034E-04 & 9.461E-02 & 9.943E-02\\
\hline
\multirow{4}{*}{PL3}
& $2^{1}$ &8.202E+10 & 1.589E+23 & -1.249E+10 & 1.765E+11\\
& $2^{0}$ &Inf & NaN & NaN & NaN\\
& $2^{-1}$&NaN & NaN & NaN & NaN\\
& $2^{-2}$&3.810E-01 & 6.677E-05 & 3.791E-01 & 3.829E-01\\
& $2^{-3}$&-9.317E-02 & 8.318E-05 & -9.533E-02 & -9.101E-02\\
& $2^{-4}$&-1.930E-02 & 9.220E-05 & -2.158E-02 & -1.703E-02\\
\hline
\multirow{4}{*}{DRI1}
& $2^{1}$ &1.360E+03 & 9.350E-02 & 1.360E+03 & 1.360E+03\\
& $2^{0}$ &3.948E+01 & 6.535E-05 & 3.947E+01 & 3.948E+01\\
& $2^{-1}$&1.525E+00 & 8.976E-05 & 1.522E+00 & 1.527E+00\\
& $2^{-2}$&-1.412E-01 & 1.012E-04 & -1.436E-01 & -1.388E-01\\
& $2^{-3}$&-3.861E-02 & 1.054E-04 & -4.105E-02 & -3.618E-02\\
& $2^{-4}$&-3.432E-03 & 1.032E-04 & -5.841E-03 & -1.022E-03\\
\hline
\multirow{4}{*}{AN3D1}
& $2^{1}$ &3.649E+01 & 8.416E-05 & 3.648E+01 & 3.649E+01\\
& $2^{0}$ &1.871E+00 & 8.809E-05 & 1.869E+00 & 1.873E+00\\
& $2^{-1}$&-4.186E-01 & 7.102E-05 & -4.206E-01 & -4.166E-01\\
& $2^{-2}$&-6.042E-02 & 6.130E-05 & -6.228E-02 & -5.857E-02\\
& $2^{-3}$&-5.103E-03 & 8.299E-05 & -7.263E-03 & -2.943E-03\\
& $2^{-4}$&3.022E-06 & 9.237E-05 & -2.276E-03 & 2.282E-03
\end{tabular}
\end{center}
\end{table}

\begin{figure}[tbp]
\includegraphics[width=6.8cm]{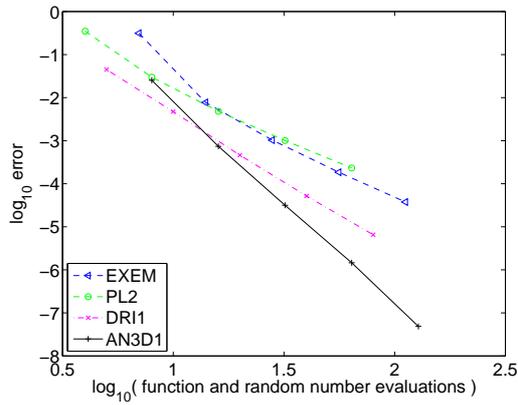}
\caption{Computational effort per simulation path versus precision for SDE \eqref{eq:ex3}} \label{aco:Bild003}
\end{figure}

As last example we consider the following linear system of SDEs with two dimensional noise
\begin{equation}\label{eq:ex3}
X(t)=\begin{pmatrix}
1\\1
\end{pmatrix}
+\int_0^t
\begin{pmatrix}
-\frac12&0\\
-\frac1{100}&-\frac34
\end{pmatrix}
X(s)~ds+\begin{pmatrix}
-\frac1{10}&\frac1{20}\\
0&\frac1{30}
\end{pmatrix}
\begin{pmatrix}
W_1(t)\\W_2(t)
\end{pmatrix},\quad t\in I=[0,2],
\end{equation}
where $\E(X_2^2(t))$ can be calculated as
\begin{equation}
    \E(X_2^2(t)) =
    \frac{37 + 31148e^{-5 t/4} - 1185 e^{-t}}{30000}.
\end{equation}
The solution value $\E(X_2^2(T))$ is approximated with step sizes $2^{1},\ldots, 2^{-3}$ at time $T=2$.
The results are presented in Table~\ref{aco:Table3} and Figure~\ref{aco:Bild003} (note that PL3 is not applicable here). Again, AN3D1 performs best.
\begin{table}
\caption{Mean errors, empirical variances and confidence intervals for SDE \eqref{eq:ex3}} \label{aco:Table3}
\begin{center}
\begin{tabular}{c|c|c|c|c|c}
   & $h$ & $\hat{\mu}$ & $\hat{\sigma}_{\mu}^2$  & $a$ & $b$ \\
   \hline
\multirow{4}{*}{EXEM}
& $2^{1}$ &-3.122E-01 & 1.279E-10 & -3.122E-01 & -3.121E-01\\
& $2^{0}$ &-7.717E-03 & 1.808E-11 & -7.718E-03 & -7.716E-03\\
& $2^{-1}$&-1.032E-03 & 2.299E-11 & -1.033E-03 & -1.030E-03\\
& $2^{-2}$&-1.848E-04 & 2.527E-11 & -1.860E-04 & -1.836E-04\\
& $2^{-3}$&-3.759E-05 & 2.876E-11 & -3.887E-05 & -3.632E-05\\
\hline
\multirow{4}{*}{PL2}
& $2^{1}$ &3.491E-01 & 9.822E-12 & 3.491E-01 & 3.491E-01\\
& $2^{0}$ &2.984E-02 & 1.040E-11 & 2.984E-02 & 2.984E-02\\
& $2^{-1}$&4.796E-03 & 8.046E-12 & 4.796E-03 & 4.797E-03\\
& $2^{-2}$&1.006E-03 & 7.765E-12 & 1.005E-03 & 1.006E-03\\
& $2^{-3}$&2.325E-04 & 6.647E-12 & 2.319E-04 & 2.331E-04\\
\hline
\multirow{4}{*}{DRI1}
& $2^{1}$ &-4.468E-02 & 3.591E-13 & -4.468E-02 & -4.468E-02\\
& $2^{0}$ &-4.712E-03 & 5.216E-12 & -4.713E-03 & -4.712E-03\\
& $2^{-1}$&-4.603E-04 & 7.853E-12 & -4.610E-04 & -4.597E-04\\
& $2^{-2}$&-5.199E-05 & 7.022E-12 & -5.262E-05 & -5.137E-05\\
& $2^{-3}$&-6.531E-06 & 6.369E-12 & -7.130E-06 & -5.933E-06\\
\hline
\multirow{4}{*}{AN3D1}
& $2^{1}$ &2.526E-02 & 9.516E-12 & 2.526E-02 & 2.526E-02\\
& $2^{0}$ &7.390E-04 & 7.190E-12 & 7.383E-04 & 7.396E-04\\
& $2^{-1}$&3.150E-05 & 5.109E-12 & 3.096E-05 & 3.204E-05\\
& $2^{-2}$&1.459E-06 & 6.996E-12 & 8.314E-07 & 2.086E-06\\
& $2^{-3}$&4.859E-08 & 6.152E-12 & -5.395E-07 & 6.367E-07
\end{tabular}
\end{center}
\end{table}
\section{Conclusion}
We have presented a general class of SRK methods for the weak approximation of SDEs with additive noise, together with the corresponding order conditions up to order three. A concrete explicit third order method has been derived, for which a numerical comparison with some well known other methods regarding its performance yielded very promising results.
In contrast to the method of Platen, it needs only two random variables and four drift evaluations per step and is also applicable to SDEs driven by a multidimensional Wiener process.
Future research may be done by constructing implicit methods with good stability properties, i.\,e.\ which are suitable for stiff problems, and by developing methods for more general noise.
\section*{Acknowledgement}
The author is grateful to Birgit Debrabant and an anonymous referee for their helpful hints which improved the presentation of the material.

\begin{thebibliography}{10}
\providecommand{\url}[1]{{#1}}
\providecommand{\urlprefix}{URL }
\expandafter\ifx\csname urlstyle\endcsname\relax
  \providecommand{\doi}[1]{DOI~\discretionary{}{}{}#1}\else
  \providecommand{\doi}{DOI~\discretionary{}{}{}\begingroup
  \urlstyle{rm}\Url}\fi

\bibitem{burrage96hso}
Burrage, K., Burrage, P.M.: High strong order explicit {R}unge--{K}utta methods
  for stochastic ordinary differential equations.
\newblock Appl. Numer. Math. \textbf{22}(1-3), 81--101 (1996).
\newblock Special issue celebrating the centenary of Runge-Kutta methods

\bibitem{burrage00oco}
Burrage, K., Burrage, P.M.: Order conditions of stochastic {R}unge--{K}utta
  methods by {$B$}-series.
\newblock SIAM J. Numer. Anal. \textbf{38}(5), 1626--1646 (electronic) (2000)

\bibitem{burrage04nmf}
Burrage, K., Burrage, P.M., Tian, T.: Numerical methods for strong solutions of
  stochastic differential equations: an overview.
\newblock Proc. R. Soc. Lond. Ser. A Math. Phys. Eng. Sci. \textbf{460}(2041),
  373--402 (2004).
\newblock Stochastic analysis with applications to mathematical finance

\bibitem{burrage99rkm}
Burrage, P.M.: {R}unge--{K}utta methods for stochastic differential equations.
\newblock Ph.D. thesis, The University of {Q}ueensland, Brisbane (1999)

\bibitem{debrabant08bao}
Debrabant, K., Kv{\ae}rn{\o}, A.: B-series analysis of stochastic
  {R}unge-{K}utta methods that use an iterative scheme to compute their
  internal stage values.
\newblock SIAM J. Numer. Anal. \textbf{47}(1), 181--203 (2008/09)

\bibitem{debrabant10ste}
Debrabant, K., Kv{\ae}rn{\o}, A.: Stochastic {T}aylor expansions: Weight
  functions of {B}-series expressed as multiple integrals.
\newblock Stoch. Anal. Appl. \textbf{28}(2), 293 -- 302 (2010).
\newblock \doi{10.1080/07362990903546504}

\bibitem{debrabant08cwa}
Debrabant, K., R\"{o}{\ss}ler, A.: Continuous weak approximation for stochastic
  differential equations.
\newblock J. Comput. Appl. Math. \textbf{214}(1), 259--273 (2008)

\bibitem{debrabant09ddi}
Debrabant, K., R\"{o}{\ss}ler, A.: Diagonally drift-implicit {R}unge-{K}utta
  methods of weak order one and two for {I}t\^o {SDE}s and stability analysis.
\newblock Appl. Numer. Math. \textbf{59}(3-4), 595--607 (2009)

\bibitem{debrabant09foe}
Debrabant, K., R\"{o}{\ss}ler, A.: Families of efficient second order
  {R}unge-{K}utta methods for the weak approximation of {I}t\^o stochastic
  differential equations.
\newblock Appl. Numer. Math. \textbf{59}(3-4), 582--594 (2009)

\bibitem{karatzas91bma}
Karatzas, I., Shreve, S.E.: Brownian motion and stochastic calculus,
  \emph{Graduate Texts in Mathematics}, vol. 113, second edn.
\newblock Springer-Verlag, New York (1991)

\bibitem{kloeden99nso}
Kloeden, P.E., Platen, E.: Numerical solution of stochastic differential
  equations, \emph{Applications of Mathematics}, vol.~21, 2 edn.
\newblock Springer-Verlag, Berlin (1999)

\bibitem{komori07mcr}
Komori, Y.: Multi-colored rooted tree analysis of the weak order conditions of
  a stochastic {R}unge-{K}utta family.
\newblock Appl. Numer. Math. \textbf{57}(2), 147--165 (2007)

\bibitem{komori07wso}
Komori, Y.: Weak second-order stochastic {R}unge--{K}utta methods for
  non-commutative stochastic differential equations.
\newblock J. Comput. Appl. Math. \textbf{206}(1), 158--173 (2007)

\bibitem{mackevicius01sow}
Mackevi{\v{c}}ius, V., Navikas, J.: Second order weak {R}unge--{K}utta type
  methods of {I}t\^o equations.
\newblock Math. Comput. Simulation \textbf{57}(1-2), 29--34 (2001)

\bibitem{milstein95nio}
Milstein, G.N.: Numerical integration of stochastic differential equations,
  \emph{Mathematics and its Applications}, vol. 313.
\newblock Kluwer Academic Publishers Group, Dordrecht (1995).
\newblock Translated and revised from the 1988 Russian original

\bibitem{roessler04ste}
R\"{o}{\ss}ler, A.: Stochastic {T}aylor expansions for the expectation of
  functionals of diffusion processes.
\newblock Stoch. Anal. Appl. \textbf{22}(6), 1553--1576 (2004)

\bibitem{roessler06rta}
R\"{o}{\ss}ler, A.: Rooted tree analysis for order conditions of stochastic
  {R}unge--{K}utta methods for the weak approximation of stochastic
  differential equations.
\newblock Stoch. Anal. Appl. \textbf{24}(1), 97--134 (2006)

\bibitem{roessler07sor}
R\"{o}{\ss}ler, A.: Second order {R}unge--{K}utta methods for {S}tratonovich
  stochastic differential equations.
\newblock BIT \textbf{47}(3), 657--680 (2007)

\bibitem{roessler09sor}
R\"{o}{\ss}ler, A.: Second order {R}unge--{K}utta methods for {I}t\^{o}
  stochastic differential equations.
\newblock SIAM J. Numer. Anal. \textbf{47}(3), 1713--1738 (electronic) (2009)

\bibitem{talay90eot}
Talay, D., Tubaro, L.: Expansion of the global error for numerical schemes
  solving stochastic differential equations.
\newblock Stoch. Anal. Appl. \textbf{8}(4), 94--120 (1990)

\bibitem{tocino02wso}
Tocino, \'{A}., Vigo-Aguiar, J.: Weak second order conditions for stochastic
  {R}unge--{K}utta methods.
\newblock SIAM J. Sci. Comput. \textbf{24}(2), 507--523 (electronic) (2002)

\end{thebibliography}

\end{document}